\documentclass[reqno]{amsart}
\usepackage{amssymb}
\usepackage[mathscr]{eucal}
\usepackage{eufrak}
\usepackage{graphicx}
\usepackage{float}
\usepackage{subfigure}
\usepackage{xspace}
\usepackage{tikz}

\newcommand{\p}{\partial}
\newcommand{\dd}{\mathrm{d}}

\newcommand{\R}{\mathbb{R}}
\newcommand{\ti}{\textsf{I}_\Omega}

\newcommand{\D}{\mathbb{\D}}

\newcommand{\LL}{\mathcal{L}}

\newcommand{\del}{{\partial}}

\newcommand{\tr}{\tilde \rho}

\newcommand{\tx}{\tilde{x}}

\newcommand{\tu}{\tilde{u}}

\newcommand{\tp}{\tilde{p}}

\theoremstyle{plain}
\newtheorem{theorem}{Theorem}[section]

\newtheorem{lemma}{Lemma}[section]

\theoremstyle{definition}
\newtheorem{definition}{Definition}[section]
\theoremstyle{remark}

\newtheorem{remark}{Remark}[section]
\numberwithin{equation}{section}

\begin{document}
\title[Measure Solutions of 1-D Piston Problem] {Measure solutions of one-dimensional piston problem for compressible Euler equations of Chaplygin gas}

\author{Aifang Qu}
\author{Hairong Yuan}

\address[A. Qu]
         {Department of Mathematics, Shanghai Normal University,
Shanghai,  200234,  China} \email{\tt afqu@shnu.edu.cn, aifangqu@163.com}

\address[H. Yuan]{School of Mathematical Sciences and Shanghai Key Laboratory of Pure Mathematics and Mathematical Practice,
East China Normal University, Shanghai
200241, China}
\email{\tt hryuan@math.ecnu.edu.cn}

\date{\today}

\begin{abstract}
We are concerned with the one-dimensional piston problem for the
compressible Euler equations of Chaplygin gas. If the piston moves
at constant subsonic speed to the uniform gas, there exists an
integral weak solution for the piston problem, consisting of a
shock separating constant states ahead of the piston. While if
the speed of the piston is sonic or  supersonic,  a singular
measure solution,  with density containing a Dirac measure supported
on the piston, shall be introduced to solve the problem. Integral weak
solution exists for the piston receding from the gas with any constant
speed, and there is no vacuum. In the extreme case as the Mach number
of the piston goes to infinity, the limiting equations and solutions
are the same as that for the polytropic gases.
\end{abstract}
\maketitle


\section{{Introduction}}\label{S:1}
This short note is devoted to understanding the following
one-dimensional isentropic compressible Euler system
\begin{equation}\label{11}
\begin{cases}
\displaystyle \del_t\rho+ \del_x(\rho u)=0,\\[8pt]
\displaystyle \del_t(\rho u)+\del_x(\rho u^2+p)=0,
\end{cases}
\end{equation}
where $t\ge0$ is time, $x\in\mathbb{R}$ is a space variable;
$\rho$,  $p$ and $u$ represent respectively  the density
of mass, scalar pressure and velocity of certain fluid flows.
These equations are derived from the law of conservation of mass
and momentum, while supposing the flow field is quite regular
(continuously differentiable). It is of great importance to know
to what extent these equations can still be mathematically meaningful
and represent some physically significant phenomena. For instance,
there is the well-known concept of admissible integral weak solutions
 (\cite{S}, see also Definition \ref{def21}), which are
Lebesgue measurable functions, usually discontinuous, satisfying
some integral relations deduced from integration-by-parts  after
multiplying \eqref{11} by compactly supported test functions, and
these integral solutions represent flow fields containing shock
waves in gas dynamics. However, previous research on Riemann problems
have shown that for Chaplygin gas, for which the state function is
\begin{equation}\label{12}
p(\rho)=-\frac{a}{\rho},
\end{equation}
where $a>0$ is a constant, concentration of mass may appear
\cite{Brenier-2005} and delta shock was introduced to be a solution.
 Such solutions are no longer Lebesgue measurable functions
on the Euclidean space. There are many works on formation, propagation
 and interactions of delta shocks now (see, for example,
\cite{CL,GLY,KW,serre,SWY,YZ}).
However, to understand the concentration phenomena and solve a general mixed initial-boundary
value problem of \eqref{11}, it is necessary to propose a general
concept of measure solutions.

To this end, in this note we study a typical mixed initial-boundary value problem of \eqref{11}, namely a piston moves at a uniform speed $V_0$ in a slim rectilinear tube with constant cross-sections; the tube fills static and homogeneous Chaplygin gas, i.e. the density is a positive number $\rho_0$ and the velocity is $0$. Recall that the local sonic speed in Chaplygin gas is given by
\begin{equation}\label{13}
 c=\frac{\sqrt{a}}{\rho}.
\end{equation}
We call
\begin{equation}\label{14}
M_0=|V_0|/c_0,\qquad (c_0=\sqrt{a}/\rho_0)
\end{equation}
the Mach number of the piston. If $M_0<1$ ($M_0=1, \ M_0>1$), the
 piston is said to be subsonic (sonic, supersonic).
We will show that if the piston moves into the gas subsonic, or
recedes from the gas, then there are integral weak solutions,
consisting of shocks and rarefaction waves respectively. However,
piecewise constant integral weak solutions do not exist if the piston
moves into the gas at sonic or supersonic speed. We then propose a
concept of measure solution to this initial-boundary value problem,
and then find a special singular measure solution for which the gas
concentrates on the piston. It is not a classical delta shock considered in
previous works. This justified rigourously the observation made in
\cite{CQ} (see Remark 1 there) on concentration of mass on the piston
when the Mach number of the piston is large. We also remark that the
classical method of constructing delta shocks
by using generalized Rankine-Hugoniot conditions does not work for
the problem considered here.

To the best of our knowledge, this is the first mathematical work on measure
 solutions of
initial-boundary value problems of Chaplygin gas with large initial
data. As a bonus, by studying measure solutions of piston problem,
which is in essence a problem on interactions of Chaplygin gases and
 physical boundaries, we see that the somewhat mysterious negative
pressure means the gas
is attracting the physical boundary, rather than pushing it away
for usually barotropic gases. See Remark \ref{rm31}.

In a previous work \cite{QYZ2}, the authors had studied the high Mach
 number limit of the piston problem for the polytropic gases.
The concept of measure solution proposed in this note is similar to
that introduced in \cite{QYZ2}. A significant difference is that, for
polytropic gas, singular measure solutions (concentration of mass on
the piston) appear only in the case that the piston moves to the gas
with Mach number $M_0=\infty$. Another difference is that if the
 piston recedes from the polytropic gas very quickly, vacuum may
appear in the tube --- while there is no vacuum for Chaplygin gas.
However,  for both the polytropic gas and Chaplygin gas, the high
Mach number limit corresponds to a vanishing pressure limit in
domains away from the piston, and the resultant limiting measure
 solutions are the same. This may be considered as one of the
reasons why Chaplygin gas might be used as an approximate model
 of polytropic gases in aerodynamics.

The hypersonic-limit problem of polytropic gas passing a
two-dimensional wedge was studied in \cite{QYZ1}, where the authors
also proposed a general concept of  measure solutions to the
two-dimensional steady non-isentropic compressible Euler equations.
In all these works, the basic idea is to relax the nonlinearity in
the Euler equations by considering all mass, momentum and pressure
to be Radon measures on the physical Euclidean spaces, and then
requiring that momentum etc. to be absolutely continuous with respect
to the measure of density, and the derivatives satisfy nonlinear
constraints  deduced from the Euler equations. This approach avoids
the confusion usually encountered by considering the product of a
Dirac measure and a discontinuous function in defining delta shocks.

We also remark that to study measure solutions of physical problems,
one shall use the partial differential equations derived
 directly from the physical principles, such as \eqref{11}. As in the
studies of shocks, nonlinear transformations of dependent variables
 might be meaningless for general integral weak solutions and measure
 solutions, hence may lead to incorrect problems for which either one
 cannot define correctly general measure solutions, or deduce
treacherous results.

The rest of the paper is organized as follows. In Section \ref{sec2}, we formulate the piston problem, and shift the coordinates system to move with the piston, thanks to the Galileo's principle of invariance for Newton's mechanics. Then we define measure solutions of the piston problem and present main theorems of this work. The theorems are proved in Section \ref{sec3} by considering each cases.

\section{The piston problem of Chaplygin gas and main results}\label{sec2}
\subsection{The piston problem of Chaplygin gas}\label{sec21}
Let the $x$-axis be the tube, and the trajectory of the piston being $x=V_0t$. Suppose the gas fills the domain $\{x<0\}$ initially,  with given constant state
\begin{equation}\label{eqid}
  U_0=(\rho,u)|_{t=0}=(\rho_0,0).
\end{equation}
Then the time-space domain to be considered is
\begin{equation}\label{eqd1}
\Omega_t=\{(t,x)~|~ x< V_0 t,\ \  t>0\}.
\end{equation}
On the piston, we impose the usual impermeable condition
\begin{equation}\label{eqbc}
\rho u(t,x)=0, \quad {\rm on}~x=V_0 t.
\end{equation}
The piston problem is to find a solution of \eqref{11}  in the domain $\Omega_t$, satisfying \eqref{eqid}\eqref{eqbc}.

We now applying the following  Galilean transformation to reformulate the piston problem:
\begin{eqnarray*}
  &&t'=t,~x'=x-V_0 t,\\
  &&\rho'(t',x')=\rho(t',x'+V_0 t'),\\
  &&u'(t',x')=u(t',x'+V_0 t')-V_0,\\
  &&p'(t',x')=p(t',x'+V_0 t'),
\end{eqnarray*}
the equations in \eqref{11} are invariant, while the domain
$\Omega_t$ is reduced to a quarter plane
\begin{equation*}\label{eqd2}
  \Omega'=\{(t',x')~|~~x'<0,~ t'>0\}.
\end{equation*}
{\it For convenience of statement,  we henceforth  consider the piston problem in this new coordinates and drop the upper index $``'"$ without confusion.} Now the domain is
\begin{equation}\label{eqd2}
  \Omega=\{(t,x)~|~~x<0,~ t>0\}.
\end{equation}
The initial condition is still \eqref{eqid},  while
the boundary condition becomes
\begin{equation}\label{eqbc2}
\rho u(t,x)=0, \quad {\rm on}~x=0.
\end{equation}

Next we carry out the following non-dimensional linear transformations
 of independent and dependent variables, which corresponds to some
similarity laws in physics:
\begin{equation}
\begin{array}{cl}\label{eqT1}
  &\begin{array}{l}\displaystyle
   \tilde t=\frac tT,\qquad \tx=\frac xL;\\
\displaystyle\tr =\frac{\rho}{\rho_0},\ \ \ \ \ \tu =\frac{u}{|V_0|}, \ \ \ \ \ \
 \tp =\frac{p}{\rho_0 V_0^2},
  \end{array}
\end{array}
\end{equation}
where $T$ and $L>0$ are constants with $L/T=V_0$.
Direct calculations show that
$\tilde{\rho} ,\tilde{u} ,\tilde{p} $
still satisfies \eqref{11}. {\it For simplicity of writing, we drop the tildes hereafter. } Then the state function of Chaplygin gas reads
\begin{equation}\label{2.10}
  {p}=-\frac1{{\rho}}\frac{1}{M_0^2},
\end{equation}
and we see that $a=1/M_0^2$. Hence the high Mach number limit $M_0\to\infty$ looks quite alike
vanishing pressure limit $p_0\to0$.

From \eqref{eqT1}, in the following  we shall take the initial data as
\begin{equation}\label{2.11}
\begin{split}
   \rho _0=1,\quad V_0= \pm1,\quad p_0=-1/M_0^2.
\end{split}
\end{equation}
The piston problem can now be rewritten as
to define and seek solutions of \eqref{11}, \eqref{2.11} and \eqref{eqbc2}
in the domain given by \eqref{eqd2}.

\subsection{Definition of measure solutions of piston problem}
The above formulation of piston problem only makes sense for classical
 solutions.   Since it turns out that the unknowns might be measures
singular to Lebesgue measure, we need
to rewrite the piston problem to be meaningful for general Radon measures. Recall that
a Radon measure $m$ on the upper plane $[0,\infty)\times\R$ could
act on the compactly supported continuous functions
$$\langle m, \phi\rangle=\int_0^\infty\int_{\R}\phi(t,x)m(\dd x\dd t)$$
where the test function $\phi\in C_0([0,\infty)\times\R)$.
One example of Radon measure is the standard Lebesgue measure $\LL^2$ on $\R^2$ . The other example is the following Dirac measure supported on a curve ({\it cf.} \cite{CL}).
\begin{definition}
\label{def31}
Let $L$ be a Lipschitz curve given by $x=x(t)$ for $t\in[0,T)$,
and $w_L(t)\in L_{\mathrm{loc}}^1(0,T)$. The Dirac measure
$w_L\delta_L$ supported on $L\subset\R^2$ with weight $w_L$ is
defined by
\begin{eqnarray}\label{eq31}
\langle w_L\delta_L, \phi\rangle=\int_0^T\phi(t, x(t))w_L(t)\sqrt{x'(t)^2+1}\,\dd t,\qquad \forall \phi\in C_0(\R^2).
\end{eqnarray}
\end{definition}
Recall that for two measures $\mu$ and $\nu$, the standard notation $\mu\ll\nu$ means  $\nu$ is nonnegative and $\mu$ is absolutely continuous with respect to $\nu$. Now we could formulate the piston problem rigorously by introducing the following definition of measure solutions.
\begin{definition}\label{def32}
For fixed  $0<M_0\leq\infty$,  let $\varrho,\ m,\  n,\ \wp$ be Radon measures on $\overline{\Omega}$, and $w_p$ a locally integrable nonnegative function on $[0,\infty)$. Then $(\varrho, u, w_p)$ is called a measure solution to the piston problem \eqref{11}\eqref{2.11}\eqref{eqbc2},  provided that
\begin{itemize}
\item[i)]
$m\ll \varrho$, $n\ll m$,  and they have the same Radon-Nikodym
derivative $u$; namely
\begin{eqnarray}\label{eq38}
u\triangleq\frac{m(\dd x\dd t)}{\varrho(\dd x\dd t)}
=\frac{n(\dd x\dd t)}{m(\dd x\dd t)};\label{eq39}
\end{eqnarray}

\item[ii)]  For any $\phi\in C_0^1(\R^2)$, there hold
\begin{eqnarray}
&&\langle \varrho, \p_t\phi\rangle+ \langle m, \p_x\phi\rangle+\int^0_{-\infty}\rho_0\phi(0,x) \dd x=0,\label{eqms1}\\
&&\langle m, \p_t\phi\rangle+ \langle n, \p_x\phi\rangle+\langle \wp, \p_x\phi\rangle-\langle
w_p\delta_{\{x=0, t\ge0\}},\phi \rangle+\int^0_{-\infty}(\rho_0u_0)\phi(0,x) \dd x=0;\nonumber\\
&&\label{eqms2}
\end{eqnarray}

\item[iii)] If $\varrho \ll \LL^2$ with derivative $\rho(t,x)$ in
 a neighborhood of $(t,x)\in[0,\infty)\times(-\infty,0]$, and
$\wp \ll \LL^2$ with derivative $p(t,x)$ there, then $\LL^2$-a.e.
there holds
    \begin{eqnarray}\label{eq312}
    p=-\frac{1}{\rho}\frac{1}{M_0^2},
    \end{eqnarray}
    and in addition, the classical entropy condition holds for
 discontinuities of functions $\rho, u$ near $(t,x)$.
\end{itemize}
\end{definition}
\begin{remark}
Physically, the weight $w_p$ is the force on the piston given by
unit volume of the gas. It is always positive when the piston moves
to the gas at supersonic speed (see \eqref{3.13}). Hence the high Mach number limit is not simply the
vanishing pressure limit, since there is an extra term
$w_p\delta_{\{x=0, t\ge0\}}$ in the limiting Euler equations (see the fourth term in \eqref{eqms2}), comparing to the
standard pressureless Euler equations (i.e., the sticky particle
 system). The requirement that $w_p$ is nonnegative shall be considered as a kind of stability condition for the measure solutions.
\end{remark}
\begin{remark}
Integral weak solutions $(\rho, u, p)$ of the piston problem  are
measure solutions, just by taking ({\it cf.} \eqref{eqdws}) $$\varrho=\rho\ti\LL^2,\ m=\rho u\ti
\LL^2,\  n=\rho u^2\ti\LL^2,\ \wp=p\ti\LL^2, \ w_p(t)\equiv0.$$
 Here $\ti$ is the
characteristic function of the set $\Omega,$ namely $\ti(t,x)=1$ if $(t,x)\in\Omega$ and $\ti(t,x)=0$ otherwise.
\end{remark}

\begin{remark}From the above definition, we may propose the following
 general formulation of Euler equations in the framework of Radon
measure solutions:
\begin{eqnarray*}\begin{cases}
\p_t\varrho+\p_xm=0,\\
\p_tm+\p_x(n+\wp)=\mathcal{F}_P,\\
\p_tm^1+\p_xn^1=0.
\end{cases}\end{eqnarray*}
Here $\mathcal{F}_P$ is a vector-valued  measure, supported on a
lower dimensional manifold $P$. It is absolutely continuous with
respect to $\varrho$. It may be used to denote the force acting on
the gas by lower dimensional surface (such as impacting considered
above, or frictions induced by physical boundary).
This formulation reduced to the classical form once the measures
in the equations share some fine regularity properties.
\end{remark}

The main results of this paper are the following two theorems:
\begin{theorem}\label{thm41}
For the piston moving towards the gases ($V_0=-1$) with sonic or
supersonic speed ($M_0\ge 1$), the problem
\eqref{11}\eqref{2.11}\eqref{eqbc2} admits a measure solution
rather than an integral
weak solution. And when $M_0=\infty$, the high Mach number limiting equations
and solutions are the same as that of polytropic gases for Euler
equations.  Integral weak solutions exist only for subsonic case
of the piston problem.
\end{theorem}

\begin{theorem}\label{thm51}
For the piston recedes from the gas ($V_0=1$), problem
\eqref{11}\eqref{2.11}\eqref{eqbc2} always has a rarefaction
wave solution, for which the wave fan degenerates to a line.
In the high Mach number limit case ($M_0=\infty$), the location
of the rarefaction wave is $x=t$. Beyond the wave, vacuum presents
 ahead of the piston. The limiting equations are the same as
that of pressureless Euler flow.
\end{theorem}

\section{Proof of main results}\label{sec3}

\subsection{Integral weak solutions}
To understand the necessity of introducing general measure solutions,
we firstly  consider integral weak solutions of the piston problem.
\begin{definition}\label{def21}
We say $(\rho,\ u)\in L^\infty([0,\infty)\times(-\infty,0])$ is
an integral weak solution to problem
\eqref{11}\eqref{2.11}\eqref{eqbc2}, if for any $\phi\in C_0^1(\R^2)$,
 there hold
\begin{equation}\label{eqdws}
 \begin{cases}\displaystyle
     \int_\Omega (\rho\p_t\phi+\rho u \p_x\phi) \dd x\dd t+
\int^0_{-\infty} \rho_0(x)\phi(0,x) \dd x=0, \\
     \begin{split} \displaystyle
      \int_\Omega (\rho u\p_t\phi+(\rho u^2+p)\p_x\phi) \dd x\dd t&
-\int_0^\infty p(t,0)\phi(t,0) \dd t\\
     &+\int^0_{-\infty} \rho_0(x)u_0(x)\phi(0,x) \dd x=0.
     \end{split}
 \end{cases}
\end{equation}
\end{definition}

\subsubsection{Shock wave solution when piston moves subsonic to
 the gas}

Noticing that problem \eqref{11}\eqref{2.11}\eqref{eqbc2} is a
Riemann problem with boundary conditions for fixed $M_0\in(0,\infty)$,
 we try to construct  self-similar  solutions  $U(t,x)=V(x/t)$.
Suppose a piecewise constant solution is of the form
\begin{equation}\label{eqshock}
U(t,x)=V(\frac{x}{t})=\begin{cases}
V_0=(1,1),& -\infty\le\frac xt<{\sigma},\\
V_1=(\rho_1,0),& {\sigma}<\frac xt\le0.
\end{cases}
\end{equation}
From \eqref{eqdws} one sees that $V_1$ and $\sigma$ shall satisfy the following Rankine-Hugoniot conditions:
\begin{equation}\label{eqrh}
  \begin{cases}
     \sigma (\rho_1-\rho_0)=\rho_1 u_1-\rho_0u_0, \\
        \sigma(\rho_1u_1-\rho_0u_0)=\rho_1 u_1^2+p_1-\rho_0u_0^2-p_0.
  \end{cases}
\end{equation}
In view of  $\rho_0=1$, $u_0=1$, $u_1=0$, it follows from
$\eqref{eqrh}_1$ that
\begin{equation}\label{eqsigma}
  \sigma=-\frac{1}{\rho_1-1}.
\end{equation}
Note that $\sigma<0$ requires that $\rho_1>1$.
Inserting it into $\eqref{eqrh}_2$ gives
\begin{equation}\label{3.6}
  p_1=p_0+1+\frac{1}{\rho_1-1}.
\end{equation}
Since $p_1=-\frac{1}{\rho_1}\frac{1}{M_0^2}=\frac{p_0}{\rho_1}$,
one has
\begin{equation}\label{eqrho1}
  (p_0+1)\rho_1^2-2p_0\rho_1+p_0=0.
\end{equation}
For $0<M_0<1$, we solve that
\begin{equation}\label{eq3.8}
   \rho_1=\frac{p_0-\sqrt{-p_0}}{p_0+1}=\frac1{1-M_0}>1.
\end{equation}
It follows from \eqref{eqsigma} that $\sigma=1-\frac{1}{M_0}$.
Therefore we proved
\begin{lemma}
  For the piston moves to the gas with Mach number $M_0<1$,
there exists an integral weak solution to
\eqref{11}\eqref{2.11}\eqref{eqbc2}. The shock wave locates at
$x=(1-\frac{1}{M_0})t$, across which the density of the gas
increases and is given by \eqref{eq3.8}.
\end{lemma}

\subsubsection{Nonexistence of shock wave solution for $M_0\geq 1$.}
For $M_0>1$, it follows from \eqref{eqrho1} and the non-negativeness of
$\rho_1$ that
\begin{equation}\label{3.7}
  \rho_1=\frac{p_0+\sqrt{-p_0}}{p_0+1}=\frac{M_0-1}{M_0^2-1}
=\frac1{M_0+1}<1.
\end{equation}
For $M_0=1$, we have $p_0=-1$, and \eqref{3.7} gives
$\rho_1=1/2$.
For these two cases $\sigma>0$ and contradicts \eqref{eqshock}.
Therefore, we conclude
\begin{lemma}
There is no piecewise constant integral weak solution to  the
piston problem of a Chaplygin gas when the piston moves to the
gas sonic or supersonic.
\end{lemma}

\subsubsection{Existence of rarefaction wave solution when the piston recedes from the gas}

Since \eqref{11} is linearly degenerate for Chaplygin gas, the rarefaction wave curves coincide with that of the shock waves in the physical plane. In view of $\rho_0=1,\ u_0=-1$, it follows from \eqref{eqsigma} that
\begin{equation}\label{3.8}
  \sigma=\frac{1}{\rho_1-1},
\end{equation}
while \eqref{3.6} is still valid now. Thus, the density $\rho_1$ is also given by \eqref{3.7}.
 It follows particularly that as $M_0\rightarrow\infty$, the density behind the rarefaction wave goes to $0$. The pressure behaves similarly, since
$$p_1=-\frac{1}{\rho_1}\frac{1}{M_0^2} =-\frac{M_0+1}{M_0^2}\rightarrow0,\quad\text{as}~M_0\rightarrow\infty.$$
The limiting location of the wave is
\begin{equation}\label{3.9}
  x=-t.
\end{equation}
Noticing that in the domain $-t\leq x\leq 0$, both the density and the pressure are zero in the limiting case. Recall the pressure ahead of the rarefaction wave is also $0$. We then conclude
\begin{lemma}
  The limiting  equations ($M_0=\infty$) to the receding piston problem of a Chaplygin gas for isentropic Euler equations are the pressureless Euler equations.
\end{lemma}

\subsection{Singular measure solutions for sonic or supersonic piston
 moves to gases} To solve the piston problem when the piston moves
sonic or supersonic to the gas, we construct a special measure
solution by supposing that
\begin{eqnarray}\label{eqms}
\varrho=\ti\LL^2+w_\rho(t)\delta_{\{x=0, t\ge0\}},\quad m=\ti\LL^2,\quad
n=\ti\LL^2,\quad
\wp=-\frac{1}{M_0^2}\ti\LL^2. 
\end{eqnarray}
Recall that $\ti$ is the characteristic function of $\Omega$. These
expressions come from the physical phenomena of infinite-thin
shock layer in
hypersonic flows past bodies and the hypersonic similarity law \cite{Anderson-2006}, and observations made in \cite[Remark 1]{CQ}.

By Definition \ref{def32}, we deduce that
\begin{equation*}
\begin{split}
    0=&\langle \varrho, \p_t\phi\rangle+ \langle m, \p_x\phi\rangle+\int^0_{-\infty}\rho_0\phi(0,x) \dd x\\
    =&\int_\Omega\p_t\phi \dd x\dd t+\int_0^\infty w_\rho(t)\p_t\phi(t,0) \dd t+\int_0^\infty\int^0_{-\infty}\p_x\phi \dd x\dd t+\int^0_{-\infty}\phi(0,x) \dd x\\
    =&\int^0_{-\infty}\phi(t,x)|_{t=0}^\infty \dd x+w_\rho(t)\phi(t,0)|_{t=0}^\infty -\int_0^\infty w_\rho'(t) \phi(t,0) \dd t+\int_0^\infty\phi(t,x)|_{x=-\infty}^0 \dd t\\&\qquad+\int^0_{-\infty}\phi(0,x) \dd x\\
    =&\int_0^\infty (1-w_\rho'(t))\phi(t,0) \dd t-w_\rho(0)\phi(0,0)\\
\end{split}
 \end{equation*}
 Due to the arbitrariness of $\phi$, we have
 \begin{equation}\label{3.11}
   \begin{cases}
    w_\rho'(t)=1&t\ge0,\\
     w_\rho(0)=0.
   \end{cases}
 \end{equation}
 It follows that
 \begin{equation}\label{3.12}
    w_\rho(t)=t.
 \end{equation}
From the momentum equation we have
\begin{equation*}
\begin{split}
    0=&\langle m, \p_t\phi\rangle+ \langle n, \p_x\phi\rangle-\frac{1}{M_0^2}\int_0^{\infty}\int_{-\infty}^0(\p_x\phi)\,\dd x\dd t-\int_{0}^{\infty}w_p(t)\phi(t,0)\,\dd t\\&\qquad+ \int^0_{-\infty}\rho_0\phi(0,x)\, \dd x\\
    =&\int_0^\infty\int^0_{-\infty}(\p_t\phi+\p_x\phi)\, \dd x\dd t-\int_0^\infty (w_p(t)+\frac{1}{M_0^2})\phi(t,0)\, \dd t+\int^0_{-\infty}\phi(0,x)\, \dd x\\
    =&\int_0^\infty(1-\frac{1}{M_0^2}-w_p(t))\phi(t,0)\,\dd t.
\end{split}
 \end{equation*}
It follows that
\begin{equation}\label{3.13}
  w_p(t)=1-\frac{1}{M_0^2}\ge0.
\end{equation}
Comparing to results in \cite{QYZ2}, in the high Mach number limiting
case $M_0=\infty$, the solution of the piston problem for the
Chaplygin gas is the same as that of polytropic gas. So we have
\begin{lemma}
There exists a measure solution for the piston problem of Chaplygin
gas when the piston moves to the gas at sonic or supersonic speed.
The limiting  equations ($M_0=\infty$) is consistent  with that
of the polytropic gases.
\end{lemma}
\begin{remark}\label{rm31}
From \eqref{3.13}, we see that for Chaplygin gas, the negative pressure means the gas is attracting the solid impermeable boundary (i.e., the piston), rather than pushing it away, comparing to polytropic gases. Particularly, for the critical case $M_0=1$, there is a balance between the impact due to the macroscopic inertia effect of a large body of particles, and the negative pressure coming from microscopical thermal motion of fluid particles, hence the piston feels no force at all (i.e., $w_p(t)\equiv0$). Without considering interactions of Chaplygin gases with physical boundaries, this attracting effect and the meaning of negative pressure would not be easy to understand.
\end{remark}
Summing up these lemmas, we finished proof of the two main theorems listed in Section \ref{sec2}.

\section*{Acknowledgements}
The research of Aifang Qu is supported by National Natural
Science Foundation of China (NNSFC) under Grant Nos. 11571357, 11871218.
Hairong Yuan is supported by NNSFC under Grant No. 11871218, and by
Science and Technology Commission of Shanghai Municipality (STCSM)
under grant No. 18dz2271000.

\end{document}